\documentclass[12pt,leqno]{amsart}
\usepackage{amsmath,amsfonts,amssymb,amsthm}

\theoremstyle{plain}
\newtheorem{theorem}{Theorem}[section]

\newtheorem{corollary}{Corollary}[theorem]
\theoremstyle{definition}

\usepackage{graphicx}
\usepackage{wrapfig}
\theoremstyle{remark}

\begin{document}
\title[]{Extension property of continuous functions in a Riemannain manifold with a pole}
\author[  A. A. Shaikh, C. K. Mondal]{Absos Ali Shaikh$^1$$^*$ and Chandan Kumar Mondal$^2$}

\begin{abstract}
The Brouwer fixed point theorem says that any continuous function from disc to itself has a fixed point. By using simple geometrical technique we have generalized the result in manifold and proved that any continuous function on the boundary of a bounded convex domain of a $2$-dimensional Riemannian manifold with a pole having at least one fixed point can be extended to the convex domain without any interior fixed point. 
\end{abstract}
\noindent\footnotetext{$^*$ Corresponding author\\
$\mathbf{2010}$\hspace{5pt}Mathematics\; Subject\; Classification: 53C22, 58C30. }
\maketitle
\section{Introduction}
The well-known classical Brouwer fixed point theorem states that any continuous function from a closed disc $D$ to itself has a fixed point. An amazing application of this theorem is the stirring of coffee. If someone stairs the coffee, then there must be some points on the top of the coffee where the points are not moving at all. But what will happen if the boundary points are fixed? One can easily stair the coffee so that the boundary points are not moving as well as there are some fixed inner points also. Naturally the question arises is there any process of stirring so that boundary points are moving but not any inner point. If the face of the coffee mug is circle, then there is a positive solution, see \cite{BG94}. But in case where the face is not circle? For any arbitrary shape the answer is still unknown. In this paper we have given a positive answer of this question in the case of bounded geodesic convex subset of a $2$-dimensional Riemannian manifold.
\par  A subset $C$ in the Euclidean plane is said to be convex if any two points can be joined by a straight line which lies in $C$. The boundary of disc $D$ and $C$ will be denoted by $S$ and $\partial C$ respectively. A subset $A$ of the Riemannian manifold $M$ is said to be geodesic convex \cite{UDR94} if for any two points, the minimal geodesic connecting them lies in $A$. A point $o$ of the Riemannian manifold is called a pole \cite{GM69} if the map $exp_o$ is a diffeomorphism from $T_oM$ to $M$. Riemannian manifold with a pole will be denoted by $(M,o)$. In $(M,o)$ any geodesic starting from $o$ is minimal \cite{IST12}. A Riemannian manifold with a pole is automatically complete.
\par The mathematical formulation of the above question can be stated as follows: Is it possible to extend any continuous function on the boundary of a convex domain with some fixed points in a continuous function of the convex domain without any interior fixed point? By using the technique of complex analysis Brown and Greene \cite{BG94} proved the following:
\begin{theorem}\label{lem1}\cite{BG94}
Given a continuous function $f:S\rightarrow S$ with at least one fixed point, there exists a map $G:D\rightarrow S$ such that $G|_S=f$ and in particular $G$ has no fixed point on $Int(D)$.
\end{theorem}
We have given an alternative geometrical proof of the result by using rotation and translation of circle. The technique developed in this proof will be used in proving the main result of the paper. By using the alternative technique of the proof of Theorem \ref{lem1}, we have given a purely geometrical proof of the following:
\begin{theorem}\label{th1}
Let $C$ be a convex subset of the Euclidean plane. Then any continuous function $f:\partial C\rightarrow\partial C$ with at least one fixed point can be extended to a continuous function $\psi:C\rightarrow \partial C$.
\end{theorem}
\begin{corollary}\label{cor1}
Let $T$ be a starshaped subset of the Euclidean plane. Then any continuous function $f:\partial T\rightarrow\partial T$ with at least one fixed point can be extended to a continuous function $\xi:T\rightarrow \partial T$. 
\end{corollary}
The main result of this paper is the following:
\begin{theorem}\label{th2}
Suppose $(M,o)$ is a $2$-dimensional Riemannian manifold with a pole $o$ and $V$ is a bounded closed convex subset of $M$ with boundary $\partial V$ such that $o\in V$. Any continuous function $f:\partial V\rightarrow\partial V$ with a fixed point can be extended to a continuous function from $V$ to $V$ without any interior fixed point.
\end{theorem}
\section{Proof of the results}

\text{\textit{Proof of Theorem \ref{lem1}:}}
Suppose that $p\in S$ is the fixed point of $f$, see Figure 1. The disc $D$ can be covered by the concentric circles $S_t$ of radius $t$ for $t\in [0,1]$, i.e., $D=\cup_{t\in[0,1]}S_t$. Now consider a continuous 

\begin{wrapfigure}{r}{0.60\textwidth}\label{fig1} 
    \centering
    \vspace{-1.2 cm}
    \includegraphics[width=0.30\textwidth]{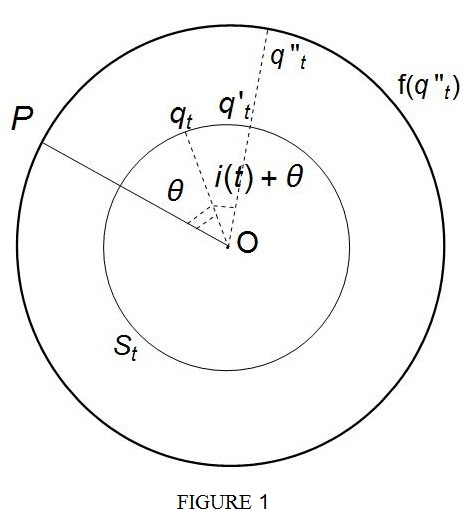}
    \vspace{-.7cm}
\end{wrapfigure}
bijection $i:[0,1]\rightarrow[0,2\pi]$ such that $i(0)=0$ and $i(1)=2\pi$. Let $q_t\in D$. Then there exists $t\in [0,1]$ such that $q_t\in S_{t}$.
Suppose $\overline{oq_t}$ makes an angle $\theta$ with $\overline{op}$. Now rotate the circle $S_{t}$ clockwise by an angle $i(t)$. The action is denoted by $\theta_{t}$. Therefore the image of $q$ is $q'_t=\theta_{t}(q_t)$ and the angle $\angle q'_top$ is $i(t)+\theta$. Now extend the line $oq'_t$ so that it intersects $S$ at the point $q''_t$. Now define the function $G:D\rightarrow D$ by
$G(q_t)=f(q''_t)\ \text{ for all }q_t\in D.$ It can easily be seen that $G(o)=p$ and for any boundary point $q\in S$, $G(q)=f(q)$. Since the function $G$ is constructed by using the composition of rotation, translation and $f$ function, $G$ is continuous. Observe that image of $G$ is the boundary set $S$. Therefore, $G$ can not have any interior fixed point. Hence, the result is proved. $\Box$
\vspace{1 cm}

\par
\begin{wrapfigure}{r}{0.50\textwidth}\label{fig2} 
    \centering
    \vspace{-.5cm}
    \includegraphics[width=0.40\textwidth]{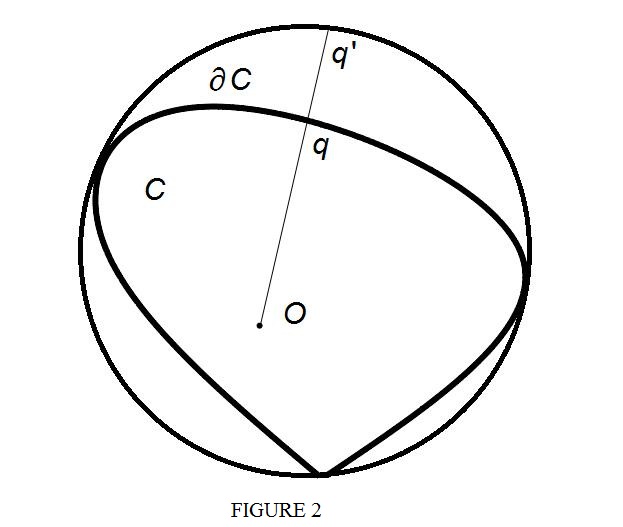}
    \vspace{-.5cm}
\end{wrapfigure}

 \text{\textit{Proof of Theorem \ref{th1}:}}
Since $C$ is a convex subset, there exists a unique circumscribed circle $S$ of $C$, see \cite{ZAM80}. Now fix a point $O$ inside $C$ such that $O\in Int(C)$, see Figure \ref{fig2}. Take a point $q\in C$, then extend the straight line from $O$ to $q$ so that it intersects $S$ at the point $q'$. Therefore, we get a continuous function $h: C\rightarrow S$, defined by $h(q)=q'$ for all $q\in C$. Again, if we restrict $h$ on $\partial C$, then we get a bijection from $\partial C$ to $S$, which will be denoted by $\partial h$. 
   Now the function $f':S\rightarrow S$ defined by
$f'(q)=\partial h\circ f\circ \partial h^{-1}(q)\ \text{ for all }q\in S.$ This map is continuous being the composition of continuous maps. Suppose $p\in\partial C$ is a fixed point of $f$ and $p'=\partial h(p)\in S$. Now $f'(p')=\partial h\circ f\circ\partial  h^{-1}(p')=\partial h\circ f(p)=\partial h(p)=p'$. This implies that $f'$ has a fixed point. Now by the Theorem \ref{lem1}, $f'$ can be extended to a continuous map $G':D\rightarrow S$ such that $G'|_S=f'$. Consider the function $\psi:C\rightarrow \partial C$ defined by
$\psi(q)=\partial h^{-1}\circ G'\circ h(q)$. If we take $q\in\partial C$, then $\psi(q)=\partial h^{-1}\circ G'\circ h(q)=\partial h^{-1}\circ G'(q')=\partial h^{-1}(f'(q'))=\partial h^{-1}\circ\partial  h\circ f\circ \partial h^{-1}(q')=f(q).$ Hence, we get a continuous function $\psi:C\rightarrow \partial C$ such that $\psi|_{\partial C}=f$.$\Box$

\vspace{1 cm}
\text{\textit{Proof of Corollary \ref{cor1}}:}
To prove it just choose the focus of $T$ as the point $O$ and run through the arguments of the Theorem \ref{th1} and Theorem \ref{lem1}. Hence the corollary is proved.$\Box$

\vspace{1 cm}

  \text{\textit{Proof of Theorem \ref{th2}:}} Since $V$ is a convex subset of $M$ containing $o$, the set $V_0=\{exp^{-1}_oq:q\in V\}\subset T_oM$ is starshpaed with focus $o$. Now consider the map $g:\partial V_0\rightarrow \partial V_0$, defined by $g(q)=exp^{-1}_o\circ f\circ exp_o(q)$ for all $q\in \partial V_0$. Suppose $p\in \partial V$ is a fixed point of $f$ and take $p'=exp^{-1}_op$.  Now
   $g(p')=exp^{-1}_o\circ f\circ exp_o(q')=exp^{-1}_o\circ f(p)=exp^{-1}_p=p'.$ So $g$ has a fixed point on $\partial V_o$. Hence in view of Corollary \ref{cor1}, we get a function $g':V_o\rightarrow \partial V_o$ such that $g'|_{\partial V_o}=g$.
\begin{wrapfigure}{r}{0.55\textwidth}\label{fig3} 
    \centering
    \vspace{0 cm}
    \includegraphics[width=0.40\textwidth]{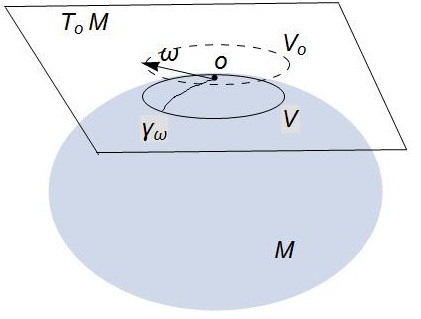}
    \vspace{ -2 cm}
\end{wrapfigure}

 Now the function $H:V\rightarrow V$ defined by
 $H(q)=exp_o\circ g'\circ exp^{-1}_o(q)$ for $q\in V$ is the extension of $f$ on $V$ and $H$ does not have any interior fixed point. Hence, the theorem is proved. $\Box$
 
$\bigskip $

$^{1}$The University of Burdwan, 

Department of Mathematics,

 Golapbag, Burdwan-713104,
 
 West Bengal, India.
 
$^1$E-mail:aask2003@yahoo.co.in, aashaikh@math.buruniv.ac.in
$\bigskip $

$^{2}$Netaji Subhas Open University,

School of Sciences,

Durgapur- 713214, West Bengal, 

$^2$E-mail:chan.alge@gmail.com

\end{document}